\documentclass[a4paper,12pt]{article}

\usepackage[italian]{babel}
\usepackage[latin1]{inputenc}

\usepackage{amsmath,amsfonts,amssymb}

\usepackage{graphicx}

\newcommand{\Q}{\ensuremath{\mathbb Q}}
\newcommand{\R}{\ensuremath{\mathbb R}}
\newcommand{\C}{\ensuremath{\mathbb C}}
\newcommand{\Z}{\ensuremath{\mathbb Z}}
\newcommand{\N}{\ensuremath{\mathbb N}}

\newcommand{\op}[4]{{}_{#4}({#1}{#2}_{#3})}

\newcommand{\s}{\varsigma}

\title{Operazione di rango zero e numeri non transitivi \\ Nulranga operacio kaj netransitivaj nombroj \\ Zeroth-rank operation and non transitive numbers}

\date{13.11.2013}

\author{Cesco Reale  \thanks{Festival Italiano di Giochi Matematici.  Sito: www.cescoreale.com/matematica  
Email: nomecognome@gmail.com, inserire "cesco" al posto di "nome" e "reale" al posto di "cognome".}  }

\begin{document}

\maketitle

\begin{abstract}

Observante la rilatojn ekzistantajn inter la elementaj operacioj kiel adicio, multipliko (iteracio de adicioj) kaj potencigo (iteracio de multiplikoj), estas difinata nova operacio (nomata   \emph{inkrementado}) kohera kun tiuj le\^goj kaj tia ke adicio rezultas iteracio de inkrementadoj. La inkrementado rezultas tre simila al la  \emph{nulacio} de Rubtsov kaj Romerio, kaj krome rezultas kohera kun  la funkcio de Ackermann. Difininte la inversan operacion de inkrementado (nomata   \emph{malinkrementado}), oni observas ke \^gi ne estas fermita en $\R$. Do, estas difinata nova aro de nombroj   (nomataj $E$,  \emph{E\^seraj nombroj}) tia ke en \^gi malinkrementado rezultas fermita. Difininte la koncepton de \emph{pse\u{u}doordigo} (analoga al ordigo, sed netransitiva), oni montras ke E\^seraj nombroj estas netransitivaj. Poste estas analizataj adicio kaj multipliko en $E$, kaj oni trovas korespondon inter $E$ kaj $\C$. Fine oni etendas  inkrementadon (kaj do pse\u{u}doordigo) al $\C$, tiel ke malinkrementado estas fermita anka\u{u} en  $\C$. 

\^Slosilvortoj en Esperanto: hiperoperacioj, inkrementado, funkcio de Ackermann, maltransitiva ordigo, netransitiva ordigo, maltransitivaj nombroj, netransitivaj nombroj, novaj nombraj aroj. \\

Observing the existing relationships between the elementary operations of addition, multiplication (iteration of additions) and exponentiation (iteration of multiplications), a new operation (named  \emph{incrementation}) is defined, consistently with these laws and such that addition turns out to be an iteration of incrementations.
Incrementation turns out to be very similar to the Rubtsov's and Romerio's \emph{zeration}, and moreover it turns out to be  consistent with Ackermann's function.
After defining the inverse operation of incrementation (named \emph{decrementation}), we observe that  $\R$ is not closed under it.
So a new set of numbers is defined (named $E$, \emph{Escherian numbers}), such that decrementation is closed on it. 
After defining the concept of \emph{pseudoorder} (analogous to the order, but not transitive), it is shown that Escherian numbers are not transitive. Then addition and multiplication on $E$ are analysed, and a correspondence between $E$ and $\C$ is found. Finally, incrementation (and hence pseudoorder) is extended to $\C$, in such a way that decrementation is closed on   $\C$ too.

English keywords: hyper-operations, incrementation, zeration, Ackermann function, intransitive order, not transitive order, intransitive numbers, non transitive numbers, not transitive numbers, new number sets. \\

Osservando le relazioni esistenti tra le operazioni elementari di addizione, moltiplicazione
(iterazione di addizioni) ed elevamento a potenza (iterazione di moltiplicazioni),
viene definita una nuova operazione (denominata \emph{incrementazione}) coerente con queste leggi
e tale che l'addizione risulti un'iterazione di incrementazioni.
L'incrementazione risulta molto simile alla  \emph{zerazione} di Rubtsov e Romerio, e risulta inoltre coerente con la funzione di Ackermann.
Definita l'operazione inversa dell'incrementazione (denominata \emph{decrementazione}),
si osserva che essa non è chiusa su $\R$.
Viene definito così un nuovo insieme numerico (denominato $E$, \emph{numeri Escheriani})
tale che su di esso la decrementazione risulti chiusa.
Definita la nozione di \emph{pseudoordinamento} (analogo all'ordinamento, ma non transitivo), si mostra che i numeri Escheriani non sono transitivi. Quindi vengono analizzate l'addizione e la moltiplicazione su $E$, e si trova una corrispondenza tra $E$ e $\C$.
Si estende infine l'incrementazione (e quindi lo pseudoordinamento) a $\C$, in maniera che la decrementazione sia chiusa anche su $\C$.

Parole chiave in italiano: iperoperazioni, incrementazione, zerazione, funzione di Ackermann, ordinamento intransitivo, ordinamento non transitivo, numeri intransitivi, numeri non transitivi, nuovi insiemi numerici.\\ 

\end{abstract}

\section{Relazioni tra le operazioni aritmetiche}

Dati $a,b,n \in \N $, definiamo

$ a \circ_1 b $ (che leggiamo $a$ composto $b$ con l'operazione di rango 1) come  $a+b$;

$ a \circ_2 b $ ($a$ composto $b$ con l'operazione di rango 2) come $a\cdot b$;

$ a \circ_3 b$ ($a$ composto $b$ con l'operazione di rango 3) come $a^b$.

Dato $a \circ b = c $, definiamo operazione inversa destra ($ \circ^{-1} $) quella che, dati $c$ e $b$,
consente di ottenere $a$ e operazione inversa sinistra ($ ^{-1} \circ $)  quella che, dati $c$ ed $a$,
consente di ottenere $b$. Vedremo che l'inversa destra va applicata a destra e l'inversa sinistra va applicata a sinistra.

Definiamo:

\[
\begin{split} 
&1.1a) \quad (a\circ)_{n+1} =  a  \, ({} \circ a)^n  =  ((a \circ a) \circ a)   ... \circ a \quad \quad  n  \enskip  volte \\
&1.2a) \quad \op{a}{\circ}{}{n+1} = (  a \circ {} ) ^n \,  a  =  a \circ (a \circ (a  \circ ... a))  \quad \quad  n  \enskip  volte \\
\end{split}
\]

ovvero $a$ composto con se stesso $n$ volte (cioè dove l'operatore $\circ$ compare $n$ volte e l'elemento $a$ compare $n+1$ volte), con composizione a destra (1.1a) o a sinistra (1.2a).

La 1.1a) è equivalente a:

\[
\begin{split}
&1.1b) \quad a  \, ({} \circ a)^n   =  (a  \, ({} \circ a)^{n-1} )  \circ  a \\
\end{split}
\]

Esempio numerico: 

$ \quad (7 \circ_1)_3 = (7+)_3= 7 (+7)^2 =  (7+7)+7 = ((7(+7)^1)+7 =  ((7( \circ_1 7)^1)  \circ_1 7 $

La 1.1b) si può riscrivere come:

\[
\begin{split}
&1.1c)  \quad a  \, ({} \circ a)^n   =  (a  \, ({} \circ a)^{n+1})   \circ ^{-1}  a \\
\end{split}
\]

Esempio numerico:

$ \quad (7 \circ_2)_2 = (7 \cdot)_2= 7 ( \cdot7)^1 =  (7 \cdot 7)+7 = ((7( \cdot 7)^2) / 7 =  ((7( \circ_2 7)^2)  \circ_2^{-1} 7 $

Dalla 1.1c) per $n=0$ si ottiene $ \quad a  \, ({} \circ a)^0   =  (a  \, ({} \circ a)^1)   \circ ^{-1}  a$, ovvero  $ \quad a   =  (a \circ a)   \circ ^{-1}  a$. Estendendo $n \in \Z$, per $n=-1$ si ottiene  $ \quad a  \, ({} \circ \, a)^{-1}   =  a \, \circ ^{-1}  a$. Per  $n=-2$ si ottiene  $ \quad a  \, ({} \circ \, a)^{-2}   =  (a  \, ({} \circ \, a)^{-1})   \circ ^{-1} \,   a$, ovvero $ \quad a  \, ({} \circ \, a)^{-2}   =  (a  \, ({} \circ^{-1} \,  a))   \circ ^{-1} \,  a =  a  \, ({} \circ^{-1} \, a)^2 $.

Quindi dalla 1.1a) si osserva che l'iterazione con composizione a destra applicata $-1$ volte equivale all'operazione inversa destra, che si applica a destra.

$ c \circ_1^{-1} b $ ($c$ composto $b$ con l'operazione inversa destra di rango 1) risulta   $c-b$;

$ c \circ_2^{-1} b $ ($c$ composto $b$ con l'operazione inversa destra di rango 2) risulta $c/b$;

$  c \circ_3^{-1} b $ ($c$ composto $b$ con l'operaz. inversa destra di rango 3) come $c \, \, \sqrt[] \,\, b $, che normalmente si scrive  $\sqrt[b]c$.

La 1.2a) invece  è equivalente a:

\[
\begin{split}
&1.2b) \quad (  a \circ {} ) ^n \,  a   =  a \circ ((  a \circ {} ) ^{n-1} \,  a) \\
\end{split}
\]


La 1.2b) si può riscrivere come:

\[
\begin{split}
&1.2c)  \quad (  a \circ {} ) ^n \,  a   =  a \, \, \, ^{-1} \!\! \circ ((  a \circ {} ) ^{n+1} \,  a) \\
\end{split}
\]

Dalla 1.2c) per $n=0$ si ottiene $ \quad  (  a \circ {} ) ^0 \,  a   =  a \, \, \, ^{-1} \!\! \circ ((  a \circ {} ) ^{1} \,  a) $, ovvero  $ \quad a   =  a \, \, \, ^{-1} \!\! \circ ( a \circ \,  a) $.
Estendendo $n \in \Z$, per $n=-1$ si ottiene  $ \quad  (  a \circ {} ) ^{-1} \,  a   =  a \, \, \, ^{-1} \!\! \circ  a $. Per  $n=-2$ si ottiene   
$ \quad  (  a \circ {} ) ^{-2} \,  a   =  a \, \, \, ^{-1} \!\! \circ ((  a \circ {} ) ^{-1} \,  a) $, ovvero $ \quad  (  a \circ {} ) ^{-2} \,  a   =  (a \, \, \, ^{-1} \! \circ) ^{2} \,  a $.

Quindi dalla 1.2a) si osserva che l'iterazione con composizione a sinistra applicata $-1$ volte equivale all'operazione inversa sinistra, che si applica a sinistra.

$ a \quad ^{-1} \!\! \circ_1 c $ ($a$ composto $c$ con l'operaz. inversa sinistra di rango 1) risulta   $c-a$;

$ a \quad ^{-1} \!\! \circ_2 c $ ($a$ composto $c$ con l'operaz. inversa sinistra di rango 2) risulta  $c/a$;

$ a \quad ^{-1}\!\! \circ_3 c $ ($a$ composto $c$ con l'operaz. inversa sinistra di rango 3) risulta  $a \, \, log \, \, c$ che normalmente si scrive    $\log_a c$.

Osserviamo che sia per il  rango 1 che per il rango 2, a causa della commutatività l'operazione inversa destra coincide con quella sinistra. Ma cambia la sintassi: le inverse sinistre dei ranghi 1 e 2 hanno sintassi inversa di quella abituale, in quanto il sottraendo e il divisore devono essere applicati a sinistra.
A partire dal rango 3 invece, a causa della non commutatività,  l'operazione inversa destra è diversa da quella sinistra.
Inoltre poiché addizione e moltiplicazione sono associative, i ranghi 2 e 3 sono generati come iterazioni di operazioni associative e dunque usare la 1.1) o la 1.2) è indifferente. Ma il rango 3 (elevamento a potenza) non è associativo, per cui a partire dal rango 4 è necessario scegliere tra le due opzioni. Per ragioni che esulano dallo scopo di questo testo, in letteratura i ranghi superiori al 3 sono definiti in genere seguendo la composizione a sinistra, per cui d'ora in avanti adotteremo questa scelta. Questo non esclude la possibilità di utilizzare sia l'inversa destra che quellla sinistra, in quanto le loro definizioni prescindono dalla scelta effettuata.

In generale, dato $m \in \N $ leggiamo $ a \circ_m b $ come $a$ composto $b$ con l'operazione di rango $m$,
dove il rango dell'operazione è definito ricorsivamente dalla sottostante relazione 2).
Possiamo constatare che le operazioni elementari di moltiplicazione ed elevamento a potenza
sono derivate dall'addizione soddisfacendo la seguente relazione:
\[
\begin{split}
&2) \quad \op{a}{\circ}{m}{n} = a \circ_{m+1} n  \\
\end{split}
\]

Esempi numerici:

$ \quad \op{7}{\circ}{1}{3} = \op{7}{+}{}{3} = 7+(7+7)= 7\cdot 3 = 7 \circ_{2} 3  $

$  \quad \op{7}{\circ}{2}{4} = \op{7}{\cdot}{}{4} = 7 \cdot  (7 \cdot (7 \cdot 7)) = 7^4 = 7 \circ_{3} 4  $


\section {Operazione di rango 0}

Ci chiediamo ora se esista e come si possa definire un'operazione di rango 0 che soddisfi le precedenti relazioni. La rappresenteremo con un operatore $ \odot $  (che chiameremo  \emph{Kis}) 
 \footnote[1]{Il simbolo $ \odot $, oltre a richiamare il simbolo 0, è già usato in fonetica per rappresentare il clic bilabiale, un suono che assomiglia a quello del bacio, e che nella variante labializzata in lingua Hadza è proprio come un bacio; perciò potremo leggere questo simbolo con il suono del bacio o con la parola \emph{Kis}, radice che in Esperanto significa appunto \emph{bacio}.}
 tale che:
\[
\begin{split}
&2.1) \quad  \op{a}{\circ}{0}{n} = \op{a}{\odot}{}{n} = a+n = a \circ_{1} n \\
\end{split}
\]

Se definiamo l'operatore $ \odot $ commutativo, l'unica definizione che soddisfa la precedente
relazione è la zerazione di Rubtsov-Romerio 
\footnote[2]{Nella versione originale (http://numbers.newmail.ru/), come operatore viene utilizzato $\circ$ e non $\odot$; poiché $\circ$ è già utilizzato in matematica per la generica legge binaria di composizione, e in questo articolo viene usato con un significato simile, associato al rango di iterazione, si è preferito usare $\odot$ anche nella 3).}:

\[
\begin{split}
3) \quad & a \odot b = a+1 \quad  a,b \in \N \quad \text{se}  \quad a > b\\
  & a \odot b = b+1 \quad a,b \in \N  \quad \text{se} \quad a < b\\
  & a \odot b = a+2 = b+2 \quad  a,b \in \N \quad \text{se} \quad  a = b\\
\end{split}
\]

Infatti dalla 2.1) abbiamo $a \odot a = a+2$ e $a \odot (a \odot a) = a \odot (a+2) = a+3 = (a+2)+1$, e quindi il risultato è il successore del maggiore, sempre tranne che per $a=b$.
Questa definizione ha infatti una particolarità: la sua estensione a $ \R $ non è continua in $a=b$.
Questo perché la relazione $ 2.1) \quad  \op{a}{\circ}{0}{n} = \op{a}{\odot}{}{n} = (a \odot)^{n-1} a = a+n = a \circ_{1} n $ ha un'incoerenza per $ n=1 $; infatti per $ n=1 $ si ha $ a = a+1 $, assurdo.

Per risolvere questo problema bisognerebbe cambiare la $2.1)$ in modo da avere  $ \op{a}{\odot}{}{1} = a $ ,
quindi$ \quad  \op{a}{\circ}{0}{n} = \op{a}{\odot}{}{n} = (a+n)-1 = (a \circ_{1} n) \quad \circ^{-1}_1 1 $. 
Allora possiamo adattare la 2) nel modo seguente:  $ \quad \op{a}{\circ}{m}{n} = (a \circ_{m+1} n) \quad  \circ^{-1}_{m+1} 1 $,  e possiamo facilmente constatare che questa nuova relazione, oltre a risolvere il problema per l'operazione di rango 0, è verificata anche per le operazioni elementari di partenza.
 Infatti: $ \op{a}{+}{}{n} = (a \cdot n)/ 1 = a \cdot n, \quad \op{a}{\cdot}{}{n} = \sqrt[1]{a^n} = a^n  $.
Sostituiamo allora la $ 2) \quad \op{a}{\circ}{m}{n} = a \circ_{m+1} n $ \quad con la 

\[
\begin{split}
& 2') \quad \op{a}{\circ}{m}{n} = (a \circ_{m+1} n) \quad  \circ^{-1}_{m+1} 1 
\end{split}
\]

 Nota: dal terzo rango in poi, le operazioni perdono la commutatività, hanno quindi
 due inverse e hanno elemento neutro (1) solo a destra; di conseguenza tra le due inverse
 solo quella destra ha elemento neutro; se nella relazione 2') avessimo considerato
 l'inversa sinistra, tale relazione non sarebbe stata valida per $m>2$; ad esempio, nel caso
 dell'elevamento a potenza (operazione di rango 3) tra radice e logaritmo abbiamo considerato
 la radice.

 La relazione da soddisfare per l'operatore Kis è ora dunque:
\[
\begin{split}
&2.1') \quad  \op{a}{\circ}{0}{n} = \op{a}{\odot}{}{n} = (a+n)-1 = (a \circ_{1} n) \quad \circ^{-1}_1 1   \\
\end{split}
\]

Se, come prima, definiamo l'operatore Kis commutativo, l'unica definizione che soddisfa la 2.1') è adesso la seguente:

\[
\begin{split}
4) \quad & a \odot b = a+1 \quad  a,b \in \N \quad \text{se}  \quad a \geq b\\
& a \odot b = b+1 \quad a,b \in \N  \quad \text{se} \quad a \leq b\\
\end{split}
\]
che è l'operazione \emph{successore del maggiore} e che chiameremo \emph{incrementazione}. 
Essa è identica alla zerazione, tranne per il fatto che è stato eliminato il problema dell'incoerenza $a=a+1$ e quindi la discontinuità per $a=b$, essa è quindi una zerazione continua (secondo la definizione propostami dallo stesso G.F. Romerio).

Questa operazione può essere banalmente estesa ai numeri reali.
L'operatore Kis è commutativo per definizione (se l'avessimo definito non commutativo, avremmo avuto l'operazione \emph{successore del primo elemento} oppure l'operazione \emph{successore del secondo elemento}, entrambe poco interessanti). Esso è non associativo, distributivo rispetto all'addizione (ma non rispetto alla moltiplicazione, né all'elevamento a potenza) e inoltre non ha elemento neutro.
Ricordiamo che in assenza di parentesi vengono effettuate prima le operazioni di rango maggiore, quindi l'addizione è prioritaria sull'incrementazione.

Non-associatività: dati ad esempio, $a<b<c-1 \in \R$ abbiamo che $(a \odot b) \odot c = b+1 \odot c = c+1 \neq a \odot (b \odot c) = a \odot c+1 = c+2$.

Distributività rispetto all'addizione: dati $a,b,c \in \R$, abbiamo $(a \odot b) + c = a+c \odot b+c $.

\section{Funzione di Ackermann}

Si può notare che questa definizione è coerente con la funzione di Ackermann.  \footnote[3]{Essa è il pi\'u semplice controesempio, dato nel 1928 da Wilhelm Ackermann, di una funzione totale (definita per ogni elemento d'ingresso), ben definita (univoca) e computabile (implementabile con un algoritmo), ma non primitiva ricorsiva (implementabile a partire da un determinato insieme di funzioni elementari di base).}
La funzione di Ackermann è una successione ricorsiva a due variabili, cos\'i definita:
\[
\begin{split}
5) \quad & A(0,n)=n+1 \quad \N \ni n \geq 0, \\
& A(m,0)=A(m-1,1) \quad \N \ni m \geq 1 \\
& A(m,n)=A[m-1,A(m,n-1)] \quad \N \ni m \geq 1, \quad \N \ni n \geq 1 \\
\end{split}
\]
Da qui si ricava: $ A(0,n)=n+1 $, $ A(1,n)=n+2 $, $ A(2,n)=2n+3 $, $ A(3,n)=2^{(n+3)}-3 $.

\begin{tabbing}

A(3,n)=$2^{(n+3)}$-3 \qquad \= \qquad 5 \= \qquad 13 \= \qquad 29 \= \qquad  61 \= \qquad 125 \= \qquad 253 \= \qquad ... \kill 

n \> 0 \> 1\> 2 \> 3 \> 4 \> 5 \> ... \\

\\

A(0,n)=n+1 \> 1 \> 2 \> 3 \> 4 \> 5 \> 6 \> ... \\

A(1,n)=n+2 \> 2 \> 3 \> 4 \> 5 \> 6 \> 7 \> ... \\

A(2,n)=2 $\cdot$ n+3 \> 3 \> 5 \> 7 \> 9 \> 11 \> 13 \> ... \\

A(3,n)=$2^{(n+3)}$-3  \> 5 \> 13 \> 29 \> 61 \> 125 \> 253 \> ... \\ 

\end{tabbing}

Definiamo ora una funzione di Ackermann modificata, che evidenzi le propriet\'a che ci interessano:
\[
\begin{split}
6) \quad & A'(0,n)=n+1 \quad \N \ni n \geq 3, \\
& A'(m,3)=A'(m-1,4) \quad \N \ni m \geq 1 \\
& A'(m,n)=A'[m-1,A'(m,n-1)] \quad \N \ni m \geq 1, \quad \N \ni n \geq 4 \\
\end{split}
\]

Da qui abbiamo: $ A'(0,n)=1+n $, $ A'(1,n)=2+n $, $ A'(2,n)=2 \cdot n $, $ A'(3,n)=2^n $.
In generale: $ A'(m,n)= 2 \circ _m n $, e quindi, in particolare, $ A'(0,n)=2 \odot n $.

\begin{tabbing}

A'(0,n)=$2 \odot n$ \qquad \= \qquad 8 \= \qquad 16 \= \qquad 32 \= \qquad  64 \= \qquad 128 \= \qquad 256 \= \qquad ... \kill

n \> 3 \> 4\> 5 \> 6 \> 7 \> 8 \> ... \\

\\

A'(0,n)=$2 \odot n$  \> 4 \> 5 \> 6 \> 7 \> 8 \> 9 \> ... \\

A'(1,n)=2+n \> 5 \> 6 \> 7 \> 8 \> 9 \> 10 \> ... \\

A'(2,n)=2 $\cdot $n \> 6 \> 8 \> 10 \> 12 \> 14 \> 16 \> ... \\

A'(3,n)=$2^n$  \> 8 \> 16 \> 32 \> 64 \> 128 \> 256 \> ... \\

\end{tabbing}

Notiamo per inciso che dalla funzione di Ackermann ricaviamo:  $ 2 \circ_m 4 = 2 \circ_{m+1} 3 $, con $ \N \ni m \geq 0$, dove le operazioni per $ m>3 $ sono definite con parentesi a destra (in accordo alla 2a): $ {}_3 (a \circ_m) = a \circ_m (a \circ_m a) $.
Inoltre, $ (2 \circ_m 2) \circ_{m+1} 1=4 $ , con $ \N \ni m \geq 0 $.

\section{Operazione inversa dell'operazione di rango 0}

Analizziamo ora l'operazione inversa di rango 0. A causa della commutatività dell'incrementazione, l'operazione inversa destra coincide con quella sinistra, quindi è unica. La chiameremo \emph{decrementazione}, e la rappresenteremo con il simbolo $ \oslash $ (che chiameremo \emph{Sik}, invertendo l'ordine delle lettere di \emph{Kis}). Per semplicità, useremo la sintassi dell'operazione inversa destra.

Possiamo constatare che:

\[
\begin{split}
7) \quad & c \oslash a = c-1 \quad a,c \in \R \quad \text{se} \quad a<c-1\\
& c \oslash a = (-\infty,c-1]  \quad  a,c \in \R \quad \text{se} \quad a=c-1\\
& c \oslash a = \nexists \quad   a,c \in \R \quad  \text{se} \quad a>c-1\\
\end{split}
\]

Nel caso $ a<c-1 $ essa coincide con l'operazione \emph{predecessore del maggiore}.
Nel caso $ a=c-1 $ essa e' indeterminata in un intervallo. Nel caso $ a>c-1 $ non esiste alcun numero reale come risultato dell'operazione.

\section{Estensioni degli insiemi numerici}

Le estensioni dei numeri naturali ($ \Z, \Q, \R, \C $) derivano tutte da una stessa esigenza.
Quella di trovare insiemi numerici su cui siano chiuse le inverse delle operazioni di addizione,
moltiplicazione ed elevamento a potenza.
Ad esempio, quando definiamo la sottrazione, operazione inversa dell'addizione, scopriamo che
essa non è sempre definita. In particolare, dato che in $ a+b=c \quad  c \geq a,b \quad
\forall a,b,c \in \N $, quando analizzo $ c-b=a $ trovo che per $ b>c \quad \nexists \quad
a \in \N $ che soddisfi la precedente relazione. Dunque la sottrazione non è chiusa su \N.
Abbiamo una serie di espressioni per ora prive di significato: 0-1, 1-2, 2-3, 0-2, 1-3, ecc.
Sappiamo per\'o che 1-0 = 2-1 = 3-2, ecc., cioè se sommo una stessa quantità a $c$ e a $b$,
il risultato dell'operazione $c-b$ non cambia. Per rispettare questa propriet\'a, raggruppiamo
le nostre espressioni in opportune classi di equivalenza, e chiamiamo -1 l'insieme delle
espressioni {0-1, 1-2, 2-3,...}, chiamiamo -2 l'insieme delle espressioni {0-2, 1-3, 2-4,...},
e cos\'i via. Definiamo in questo modo \Z, un ampliamento di \N \quad su cui è chiusa
anche la sottrazione.
Analogamente dalla divisione otteniamo \Q ~ e da radice e logaritmo otteniamo \R ~ e \C ~
(in realt\'a in questa maniera non si ottengono i trascendenti, ma questo adesso non ci interessa).

\section{Definizione di un nuovo insieme numerico}

Torniamo ora al Sik. Vogliamo cercare un'opportuna estensione che renda chiuso il Sik sul nuovo insieme, e che rispetti il significato del Kis.
Dobbiamo definire un $ b $ tale che $ a \odot b = c $ quando $ a>c-1 $.  Poiché Kis è un operatore \emph{successore}, per far s\'i che $ c $ sia, in un senso nuovo, il successore di uno degli operandi, possiamo definire $ b $ come $ \s (c-1) $, dove $ \s $ (stigma) è un operatore unario scelto per rappresentare questi nuovi numeri, tale che $ a \odot \s(b) $ sia ora il \emph{successore del minore} tra $ a $ e $ b $.
 \footnote[4]{Nota: lo stigma è una lettera, pronunciata [st], che nel greco medievale e fino al XIX secolo veniva usata come legatura di sigma e tau e come cifra per indicare il 6; essa è quasi identica al sigma usato in fine di parola.}
Chiameremo questi numeri \emph{stigmareali} ($\s \R$).
Con questa nuova definizione abbiamo:

\[
\begin{split}
8) \quad & c \oslash a = c-1 \quad a,c \in \R \quad  \text{se} \quad a<c-1\\
& c \oslash a = (-\infty,c-1] \cup [\s(c-1), \s \infty) \quad  a,c \in \R \quad \text{se} \quad a=c-1\\
& c \oslash a = \s (c-1) \quad  a,c \in \R \quad \text{se} \quad a>c-1\\
\end{split}
\]

Questa definizione equivale a fare il seguente raggruppamento in classi di equivalenza delle espressioni $ c \oslash a $ con $ a>c-1 $: $ c \oslash (c+d) = c \oslash (c+f) \quad   d,f \in \R $ e $ d,f>-1 $, in analogia al caso $ a<c-1 $, dove $ c \oslash (c+d) = c \oslash (c+f) \quad   d,f \in \R $ e $ d,f<-1 $.
Per rendere la definizione univoca, nel caso $ a=c-1 $ all'interno dell'insieme
$ (-\infty,c-1]  \cup [\s(c-1), \s \infty) $ scegliamo come valore principale $ c-1 $,
analogamente a come si sceglie $ \arcsin 0 = 0 $ all'interno dell'insieme $ k \pi $,
con $ k \in \Z $. In seguito ripeteremo varie volte questa operazione: restringeremo
un insieme di risultati ad un solo risultato (che chiameremo appunto \emph{valore principale})
per rendere la funzione univoca.
Abbiamo così:

\[
\begin{split}
8.1) \quad & c \oslash a = c-1 \quad a,c \in \R \quad  \text{se} \quad a \leq c-1\\
& c \oslash a = \s (c-1) \quad  a,c \in \R \quad  \text{se} \quad a>c-1\\
\end{split}
\]

Mantendendo la commutatività, ridefiniamo ora piú precisamente l'incrementazione nella maniera seguente:

\[
\begin{split}
9) \quad & a \odot b = b  \odot a = \{a+1, \s a + 1 \} \quad a,b \in \R \quad  \text{se} \quad a \geq b\\
& a \odot b =  b  \odot a = \{b+1, \s b + 1\} \quad a,b \in \R \quad  \text{se} \quad a \leq b\\
& a \odot \s b = \s b  \odot a = \{\s b+1, b + 1\} \quad a,b \in \R \quad  \text{se} \quad a \geq b\\
& a \odot \s b = \s b  \odot a = \{a+1, \s a + 1\} \quad a,b \in \R \quad  \text{se} \quad a \leq b\\
\end{split}
\]
Ovverosia, definiamo come risultato dell'incrementazione sempre due valori, uno è il successore di un reale e l'altro è il successore di uno stigmareale.
Il valore principale è il primo dei due valori mostrati nella 9), ed è il successore dell'operando scelto.
Possiamo dunque asserire che aggiungere un operatore stigma inverte la scelta dell'operando di cui fare il successore.
In questo modo possiamo definire anche :
\[
\begin{split}
9.1) \quad & \s a \odot \s b = \s b  \odot  \s a = \{ \s a+1, a + 1 \} \quad a,b \in \R \quad  \text{se} \quad a \geq b\\
& \s a \odot \s b = \s b  \odot  \s a  = \{\s b+1, b + 1\} \quad a,b \in \R \quad  \text{se} \quad a \leq b\\
\end{split}
\]

\[
\begin{split}
9.2) \quad & \s(\s a)=a \\
\end{split}
\]

\section{Pseudoordinamento e numeri Escheriani}

Per rendere il Sik chiuso sul nuovo insieme numerico bisogna ancora definire
$ \s c \oslash a $ e $ \s c \oslash \s a $, oltre che $ \odot$ e $\oslash $ in $ \R \cup \s \R $.
Prima di ciò, definiamo la relazione  $ \succ ,\asymp,\prec $ che chiameremo \emph{pseudomaggiore}, \emph{pseudouguale}, \emph{pseudominore} e che tra reali coincide con \emph{maggiore}, \emph{uguale}, \emph{minore}.
Estendiamo la definizione dell'operatore Kis da \emph{successore del maggiore} (in $\R$)
a \emph{successore dello pseudomaggiore} (in $\R \cup \s \R$); dato $ a<b $ ($a,b \in \R$), abbiamo $ \s a \odot b = \s a+1 $,
ed essendo $\s a+1$ successore di $\s a$ definiamo $ \s a  \succ b $.
Invece, $ a>b \Rightarrow \s a \prec b $.
Inoltre $a<b \Rightarrow \s a \prec \s b$, in quanto $ \s a \odot \s b = \s b+1 $, dunque lo pseudomaggiore è $ \s b$. Infine, definiamo $a \asymp \s a$ (vedremo in seguito che questa scelta è obbligata sia per mantenere la validità di alcune formule, sia per estendere lo pseudoordinamento a $\C$).
Quindi, dato $a<b$, troviamo che $a \prec b \prec \s a \prec \s b \prec a$.
Si può notare che esistono infinite n-uple di numeri in $ \R \cup \s\R $, per le quali
$ p_1 \prec p_2 \prec ... \prec p_n \prec p_1  $. Questo ricorda alcune litografie di Escher, come \emph{Cascata} e \emph{Salita e discesa}, che sfruttano l'illusione ottica del triangolo di Reutersvärd-Penrose.
Perciò ho chiamato l'insieme $ \R \cup \s \R$  \emph{numeri Escheriani} ($E$).
Questo pseudoordinamento è dunque riflessivo, antisimmetrico e non transitivo.

Possiamo ora riscrivere la definizione di incrementazione in maniera compatta, sfruttando la nozione di pseudoordinamento:
\[
\begin{split}
10) \quad & x \odot y= \{x+1, \s x +1 \} \quad  x,y \in E  \quad \text{se}  \quad x \succeq y\\
& x \odot y =\{y+1, \s y +1 \}  \quad  x,y \in E  \quad \text{se} \quad x \prec y\\
\end{split}
\]
dove dei due risultati il valore principale è il successore dello pseudomaggiore (rispettivamente $x+1$ e $y+1$ nei due casi della 10), dunque se lo pseudomaggiore è reale il valore principale è il successore di un reale, e se lo pseudomaggiore è stigmareale il valore principale è il successore di uno stigmareale.

\section{Decrementazione in $E$ }

Torniamo ora al Sik. Analizziamo $ c \oslash \s a $.
\[
\begin{split}
11) \quad & c \oslash \s a = \s (c-1) \quad   a,c \in \R \quad  \text{se} \quad a<c-1\\
& c \oslash \s a = (-\infty,c-1]  \cup [\s(c-1), \s \infty) \quad  a,c \in \R \quad \text{se} \quad a=c-1\\
& c \oslash \s a = c-1 \quad   a,c \in \R \quad  \text{se} \quad a>c-1\\
\end{split}
\]
Anche in questo caso la rendiamo univoca, ponendo:
\[
\begin{split}
11.1) \quad & c \oslash \s a = \s (c-1) \quad    a,c \in \R \quad \text{se} \quad a \leq c-1\\
& c \oslash \s a = c-1 \quad   a,c \in \R \quad  \text{se} \quad a > c-1\\
\end{split}
\]

\section{Addizione tra Escheriani}

Cerchiamo ora di capire come funziona l'addizione in $E$.
Volendo estendere la distributività dell'incrementazione anche ai numeri stigma, dati ad esempio $a<b<c \in \R$ (negli altri casi si ottiene comunque la 12), troviamo $(a \odot \s b) + c = a+c \odot \s b+c $.
Ora, il primo membro è uguale a $a+1+c$, e nel secondo membro compare $\s b+c$, ma non sappiamo come si sommano reali e stigmareali; se poniamo $ \s b + c = \s (b+c)$, il secondo membro diventa uguale al primo.
Inoltre troviamo: $( a \odot \s b) + \s c = a+ \s c \odot \s b+ \s c $, che, in base a quanto appena detto, e mantenendo l'addizione commutativa, diventa: $  \s (a+1+c) = \s (a + c) \odot \s b+ \s c$.
Affinché quest'espressione sia valida, deve risultare  $\s b + \s c = b+c$.
Riassumendo, definiamo l'addizione in $\R \cup \s \R$ in questa maniera:

\[
\begin{split}
12) \quad & a + \s b = \s b + a =  \s (a+b) \quad    a,b \in \R \\
& \s a + \s b = \s b + \s a = a+b \quad    a,b \in \R \\
\end{split}
\]

\section{Chiusura della decrementazione in $E$ }

Tornando al problema di rendere chiuso il Sik su tutto $E$, possiamo farlo in una maniera molto semplice: basta riprendere la 8.1) e la 11.1) e considerarle valide anche quando il primo operando è stigmareale.
Analogamente all'operatore unario \emph{modulo} che elimina il segno -, definiamo ora l'operatore unario \emph{stigmamodulo} di $x$ ($\nmid x \nmid$) che elimina il "segno stigma", cioè l'eventuale presenza dell'operatore Stigma: $\nmid x \nmid ~ =  ~ \nmid \s x \nmid  ~ = x$.
Ridefiniamo a questo punto la decrementazione in maniera piú compatta:

\[
\begin{split}
13) \quad & z \oslash x =  ~ \nmid x+z-1 \nmid  ~ - x     \quad   x,z \in E \quad \text{se} \quad \nmid x \nmid  ~ \leq  ~ \nmid z-1 \nmid\\
          & z \oslash x =  ~ \nmid x+z-1 \nmid  ~ - \s x  \quad   x,z \in E \quad  \text{se} \quad \nmid x \nmid  ~ > ~  \nmid z-1 \nmid\\
\end{split}
\]

\section{Parallelismi tra formule}

Osserviamo che la 12) risulta in perfetta analogia con $a \cdot (-b) = (-b) \cdot a = -(a \cdot b)$ e $-a \cdot (-b) = -b \cdot (-a) = a \cdot b$.
In pratica basta alzare ogni operatore di un rango: $\odot$ diventa $+$, $+$ diventa $\cdot$, e l'operatore unario $\s$ diventa l'operatore unario $-$ (quello con cui distinguiamo $+a$ e $-a$, e che è diverso dall'operatore binario $-$ che ci consente di fare $c-a=b$).

Possiamo constatare che in molti casi alzando o abbassando il rango $m$ degli operatori di una formula, essa resta valida.
Abbiamo visto un esempio sull'addizione, analizzeremo ora altri esempi, e questo ci sarà utile nel prossimo paragrafo per studiare la moltiplicazione in E. Consideriamo $k,m \in \N, ~ a,b,c,d \in \R$.\\

Somma di $k$ termini consecutivi di una progressione aritmetica e prodotto di $k$ termini consecutivi di una progressione geometrica:

\[
\begin{split}
14.0) \quad & a_1  ~ \circ_m  ~ a_2  ~ \circ_m ... \circ_m  ~ a_k = (a_1  ~ \circ_m  ~ a_k) \circ_{m+1}  ~ k ~ \circ^{-1}_{m+1}  ~ 2 \\
14.1) \quad & a_1 + a_2 + ... + a_k = (a_1 +a_k) \cdot k/2 \\
14.2) \quad & a_1 \cdot a_2 \cdot ... \cdot a_k = \sqrt[2] {(a_1 \cdot a_k) ^ k} \\
\end{split}
\]
Alzando o abbassando ulteriormente di un rango la formula perde di validità.

Proprietà delle potenze e distributività:

\[
\begin{split}
15.0) \quad & a ~ \circ_{m+1}  ~ b  ~ \circ_m ~ c ~ \circ_{m+1} ~  b = (a  ~ \circ_{m}  ~ c)  ~ \circ_{m+1}  ~ b \\
15.1) \quad & a^b \cdot c^b = (a \cdot c) ^ b   \\
15.2) \quad & a \cdot b + c \cdot b = (a+c) \cdot b\\
15.3) \quad & a+b \odot c+b = (a \odot c) + b\\
\end{split}
\]

Proprietà delle potenze e commutatività:

\[
\begin{split}
16.0) \quad & (b  ~ \circ_m  ~ c) \circ^{-1}_m  ~ a = (b  \circ^{-1}_m  ~ a) \circ_m  ~ c \\
16.1) \quad & \sqrt[a] {b ^ c} = (\sqrt[a] {b}) ^ c  \\
16.2) \quad & (b \cdot c)/a =  b/a  \cdot c \\
16.3) \quad & (b + c) - a = (b-a) + c \\
\end{split}
\]
Vedremo nella 19) che la proprietà analoga partendo dal logaritmo (anziché dalla radice) non è generalizzabile.

In alcuni casi lo slittamento non è perfetto, in quanto un operatore non cambia, come negli esempi seguenti. Per questa ragione in questi esempi non è possibile estrapolare una formula valida per ranghi diversi. Possiamo però osservare che tra le formule che è possibile variare di rango (come nelle 14-21), gli operatori che restano invariati sono quelli il cui risultato è un operando destro di operazoni (dirette o inverse) di rango superiore al secondo (come nelle 17, 20, 21) oppure quelli per cui almeno un operando è il risultato di un operazione inversa sinistra di rango superiore al secondo (come nelle 18, 19). \\
Utilizzando la notazione polacca inversa, e considerando che gli operatori commutativi rendono possibili diverse scritture per la stessa formula, possiamo anche esprimere quanto sopra nel modo seguente: tra le formule che è possibile variare di rango (come nelle 14-21), gli operatori che restano invariati sono quelli che in notazione polacca inversa precedono immediatamente un operatore (diretto o inverso) di rango superiore al secondo (come nelle 17, 20, 21) oppure, in almeno una delle possibili scritture, seguono immediatamente un operatore inverso sinistro di rango superiore al secondo (come nelle 18, 19). \\

Prodotto binomiale:

\[
\begin{split}
17.1) \quad & (a \cdot b)^{c+d}=a^c \cdot a^d \cdot b^c \cdot b^d \\
17.2) \quad & (a+b) \cdot (c+d) = a \cdot c + a \cdot d + b \cdot c + b \cdot d \\
17.3) \quad & (a \odot b) + (c \odot d) = (a+c \odot a+d) \odot  (b+c \odot b+d)\\
\end{split}
\]
Si noti che, data la non-associatività dell'incrementazione, nella 17.3 si rende necessaria l'introduzione di parentesi. Un inserimento differente farebbe decadere la validità della formula.

Logaritmi e distributività (nel secondo membro delle prime due righe il + resta invariato):

\[
\begin{split}
18.1) \quad & log_a (b \cdot c) = log_a b + log_a c \\
18.2) \quad & (b+c)/a = b/a + c/a \\
18.3) \quad & (b \odot c) - a = b-a \odot c-a\\
\end{split}
\]

Logaritmi e associatività-commutatività (nel secondo membro delle prime due righe il $\cdot$ resta invariato):

\[
\begin{split}
19.1) \quad & log_a (b ^ c) = c \cdot log_a b   \\
19.2) \quad & (b \cdot c)/a = c \cdot b/a \\
19.3) \quad & (b + c) - a = c + (b-a)\\
\end{split}
\]
Abbassando ulteriormente di un rango la formula perde di validità.

Somma di frazioni (nel secondo membro delle prime due righe $b \cdot d$ resta invariato):

\[
\begin{split}
20.1) \quad & \sqrt[b]a \cdot \sqrt [d]c = \sqrt [b \cdot d] {a^d \cdot c^b}   \\
20.2) \quad & a/b + c/d = (a \cdot d + b \cdot c) / (b \cdot d)\\
20.3) \quad & a-b \odot c-d = (a+d \odot b+c) - (b+d)\\
\end{split}
\]

Per la prossima formula devo anticipare che  indico l'operazione $ \circ_4$, la tetrazione, con il simbolo $\nearrow $, visto che considero poco adatte entrambe le notazioni normalmente usate (la notazione $a \uparrow \uparrow b$ e la notazione $^b a$).

\[
\begin{split}
21.0) \quad & a ~ \circ_{m+2}  ~ b  ~ \circ_m  ~ a =  a  ~ \circ_{m+1} ~ (a  ~ \circ_{m+2}  ~ (b-1) + a ~~ ^{-1} \!\! \circ_{m+1} ~  a) \\
21.1) \quad & a \nearrow b \cdot a =  a^{(a \nearrow (b-1) + log_a a)} \\
21.2) \quad & a^b+a = a \cdot (a^{b-1}+a/a)\\
21.3) \quad & a \cdot b \odot a = a + (a \cdot (b-1) \odot a-a)\\
\end{split}
\]

Occorre notare che in tutti questi casi le formule restano valide anche con le rispettive operazioni inverse, ivi compresa la decrementazione, anche quando compaiono numeri stigmareali. Ad esempio se nella 20.3)  sostituiamo il Kis col Sik, nell'ipotesi $c-d>a-b-1$ otteniamo da entrambi i lati $\s (a-b-1)$.
Inoltre, un'analisi dei vari casi della 21.3) mostra che essa resta valida solo se valgono le condizioni finora viste su $E$: ad esempio, solo se l'incrementazione ha due risultati, solo se il valore principale è il successore dello pseudomaggiore, solo se il Sik rispetta la 13), solo se l'addizione rispetta la 12).

Bisogna inoltre evidenziare un caso particolare: analizziamo la 21.3) per $a= \s 0,  b=0$, e la 17.3) quando $b= \s a$ e $ c \odot d$ è stigmareale (o, equivalentemente,  la distributività nel caso $(a \odot \s a ) + \s c $). In questi casi, per mantenere la validità delle formule, dobbiamo scegliere attentamente quale sia il valore principale di $ a \odot \s a $ tra $a+1$ e $\s a +1$: notiamo che $ a \odot \s a $ non può essere uguale a $ \s a \odot  a $, ma deve esserlo in stigmamodulo, quindi deve essere  $ \s a \odot  a = \s ( a \odot \s a )$. In altre parole, quando i due operandi sono pseudouguali l'operazione   \emph{successore dello pseudomaggiore} deve diventare o \emph{successore del primo elemento} o \emph{successore del secondo elemento}. Noi l'abbiamo definita nel primo modo (nella 10). Questa è la prima ragione per cui  $\s a$ non può essere né pseudomaggiore né pseudominore di $a$. Nell'estensione a $\C$ vedremo una seconda ragione.

\section{Moltiplicazione tra Escheriani}

Cerchiamo ora di capire come funziona la moltiplicazione con i numeri Escheriani.
Dalla 12) ricaviamo che da $\s a \cdot b = \op{\s a}{+}{}{b}$ otteniamo un reale se $b$ è pari e uno stigmareale se $b$ è dispari, ma se $b$ non è intero cosa succede ?
Il problema è analogo al caso $(-a)^b$ con $a$ positivo: se $b$ è pari il risultato è positivo, se $b$ è dispari il risultato è negativo, se $b$ non è intero bisogna ricorrere ai numeri complessi.
Proviamo allora ad analizzare la 21.3), estendendone la validità a E.

Studiamo il caso $ \s a \cdot b \odot \s a =  \s a + ( \s a \cdot (b-1) \odot  \s a - \s a)$.
Si constata che questa formula è coerente con quanto discende dalla 12). Infatti dato ad esempio $a(b-1) \geq 0$ e $b$ pari, il primo membro risulta $ab \odot \s a = \s a+1$ e il secondo risulta $\s a + (\s (a \cdot (b-1)) \odot 0 )= \s a + 1 $. Anche negli altri casi la formula resta valida.


Per analizzare la moltiplicazione con razionali, cambiamo leggermente la 21.2) nel modo seguente: dati $b,m,n \in \Z$, $a^b+a^{m/n} = a^{m/n} \cdot (a^{b-m/n}+a^{m/n}/a^{m/n})$; abbassando di un rango, e sostituendo $a$ con $\s a$ la 21.3) diventa: $ \s a \cdot b \odot a \cdot m/n = \s a \cdot m/n + (\s a \cdot (b-m/n) \odot a \cdot m/n -a \cdot m/n )$.
Detto $b-m/n$ = $k/n$, possiamo verificare che, affinché la formula resti valida il prodotto $\s a \cdot k/n$ risulta: stigmareale ($\s (ak/n)$) per $k$ ed $n$ dispari, reale ($ak/n$) per $k$ pari ed $n$ dispari, non definito per $k$ dispari ed $n$ pari. Infatti in questo terzo caso, arriviamo a una contraddizione sia nell'ipotesi $\s a \cdot k/n = \s (ak/n)$, sia nell'ipotesi $\s a \cdot k/n = ak/n$; questo è coerente col fatto che addizionando $n$ volte (con $n$ pari) uno stesso numero (reale o stigmareale) il risultato è reale. Il tutto è ancora in analogia con $(-a)^b$.

Ora studiamo il caso $ \s a \cdot  \s b \odot \s a =  \s a + ( \s a \cdot ( \s b-1) \odot  \s a - \s a)$.
Analizziamo ad esempio il caso $a(b-1) \geq 0$ (tutti gli altri casi portano agli stessi risultati).  Se al primo membro supponiamo $\s a \cdot \s b = ab$, allora al primo membro otteniamo $\s a+1$ e affinché il secondo membro sia uguale al primo deve essere $\s a \cdot \s (b-1)  \odot 0 = 1$, e quindi $\s a \cdot \s (b-1) = \s (a \cdot (b-1))$.
Se invece al primo membro supponiamo  $\s a \cdot \s b = \s ab$ allora al primo membro otteniamo $\s ab+1$ e affinché il secondo membro sia uguale al primo deve essere $\s a \cdot \s (b-1)  \odot 0 = a(b-1)+1$, e quindi $\s a \cdot \s (b-1) = a \cdot (b-1)$.

\[
\begin{split}
22.1) \quad & \s a \cdot \s b = ab \Leftrightarrow \s a \cdot \s (b-1) = \s (a \cdot (b-1)) \\
22.2) \quad & \s a \cdot \s b = \s ab \Leftrightarrow \s a \cdot \s (b-1) = a \cdot (b-1)\\
\end{split}
\]

Le 22.1) e 22.2) sono quindi equivalenti: considerando il prodotto $ \s x \cdot \s y$, un cambio di parità di $x+y$ determina un cambiamento del risultato tra $xy$ e $\s(xy)$. Non si riesce però a ricavare se, dato $a+b$ pari,  $\s a \cdot \s b$ risulti reale o stigmareale.

Viste le analogie con il caso $(-a)^b$, analizziamo la formula generale per l'elevamento a potenza di numeri negativi e proviamo ad abbassarla di un rango (come visto alla fine del paragrafo sui parallelismi, l'operatore unario $-$ diventa l'operatore unario $\s$); dati $a \in \R^+, b \in \R, k \in \Z$ otteniamo:

\[
\begin{split}
23.1) \quad & (-a)^b = a^b \cdot e^{i \pi b(2k+1) ~ mod ~ i2 \pi} \\
23.2) \quad & \s a \cdot b = a \cdot b + e \cdot (i \pi b(2k+1) ~ mod ~ i2 \pi) = a \cdot b + i \pi e b(2k+1) ~ mod ~ i2 \pi e \\
\end{split}
\]

La 23.2) è valida per $a \in \R$, non solo in $\R^+$.
L'asse immaginario è dunque ristretto a $[0, i 2 \pi e)$, e ogni valore $iy$ esterno a questo settore è equivalente al valore principale $iy ~mod ~i2 \pi e$.
Questo è in perfetta analogia con quanto avviene per i logaritmi dei numeri complessi, dove la parte immaginaria del risultato è definita su $i\R / i2\pi$.
Abbiamo dunque che  $ \s a = a \pm i \pi e$. In particolare, troviamo l'identità di Eulero per gli stigmareali:

\[
\begin{split}
24.1) \quad & e^{\pm i \pi} =-1 \\
24.2) \quad & \pm ei \pi = \s 0 \\
\end{split}
\]

Possiamo constatare che la 23.2) rispetta le proprietà viste finora degli stigmareali. Infatti: \\
$ \quad \quad a+ \s b = a + b \pm i\pi e= \s (a+b) \quad$ e $ \quad \quad \s a + \s b = a + b \pm i2\pi e= a+b$.
Inoltre, coerentemente con le 22): $\s(\s a)= \s (a \pm i \pi e) = a \pm i2 \pi e= a$.
Analizzando la 23.2) troviamo che, in analogia alle $n$ radici n-esime di un numero complesso, anche qui ogni stigmareale, e piú in generale ogni numero complesso, ha $n$ quozienti n-esimi: $\s a/n$ corrisponde a $n$ diversi valori con parte reale $a/n$ e parti immaginarie $ (2k+1) \pi e /n$.

Per rendere univoca la 23.2) definiamo come valore principale quello per $k=0$. D'ora in avanti per semplicità considereremo solo il valore principale.
A questo punto la moltiplicazione in $E$ risulta cosí definita:
\[
\begin{split}
25.1) \quad & \s a \cdot b = b \cdot \s a = a \cdot b + i \pi e b ~ mod ~ i2 \pi e \\
25.2) \quad & \s a \cdot \s b = \s b \cdot \s a = (a + i \pi e) \cdot (b + i \pi e) = ab - (\pi e) ^2 + i \pi e (a+b) ~ mod ~ i2 \pi e\\
\end{split}
\]
Dunque adottando la 23.2) risulta che $\s a \cdot \s b $ è reale se $a+b$ è pari, e stigmareale se $a+b$ è dispari.

Osserviamo che un'idea simile alla 23.2) è stata analizzata da Rubtsov e Romerio in \emph{Progress report on hyper-operations: zeration} (2007), anche se con esiti diversi.
L'idea è quella di considerare la formula del logaritmo dell'inverso, e abbassare il rango di alcuni operatori (ma non del logaritmo):
\[
\begin{split}
26.1) \quad & ln (x^{-1}) = ln (1/x) = -ln (x)   \\
26.2) \quad & ln(x \cdot -1) = ln (0-x) = \s ~ ln (x)\\
\end{split}
\]
Dalla 26.2) deriverebbe $ \s a = a + i \pi (2k+1) ~ mod ~ i2 \pi$.

Per il resto, i numeri Delta definiti da Rubtsov e Romerio presentano molte similitudini con gli Stigmareali, pur mancando, tra le altre cose, la nozione di non transitività e l'identificazione con i numeri complessi.

\section{Pseudoordinamento, incrementazione e \\ decrementazione in $\C$}

Possiamo ora provare ad estendere lo pseudoordinamento ai complessi.
Dati $ a \in \R^+, \quad d,f \in \C,  \quad  d=d_r+id_i$ e $f=f_r+if_i$ con $d_r, f_r \in \R,  \quad  d_i, f_i \in \R /2\pi e $, le condizioni da soddisfare sono:
\[
\begin{split}
27.1) \quad & f+d \succ f \Leftrightarrow f-d \prec f\\
27.2) \quad & f+a \succ f \\
27.3)  \quad & d \succ f \Leftrightarrow d+i\pi e \prec f\\
\end{split}
\]
La 27.1) serve a mantenere lo pseudoordinamento antisimmetrico. La 27.2) è necessaria per rispettare  l'ordinamento su \R.
 La 27.3) è necessaria per rispettare lo pseudoordinamento su $E$. 

La più semplice e naturale estensione di pseudoordinamento su $\C$ coerente con le proprietà viste finora è la seguente.
Si sposti opportunamente la fascia di altezza $i 2\pi e$ (chiusa in alto e aperta in basso) in modo da centrarla su uno dei due numeri $d,f$ da comparare (ad esempio $f$). Con questi valori di $d$ e $f$ nella fascia scelta possiamo stabilire che:

\[
\begin{split}
28) \quad &  d \asymp f \Leftrightarrow ~  d-f \in \{0, i \pi e \} \\
      \quad &  d \succ f \Leftrightarrow d_i-f_i<0  \qquad \qquad \qquad \qquad   \quad   \text{se}  \quad  d_r=f_r   \\
      \quad &  d \succ f \Leftrightarrow  -\pi e/2  \leq d_i-f_i < \pi e/2   \quad \quad \quad  \text{se}  \quad  d_r > f_r   \\
      \quad &  d \prec f \Leftrightarrow  -\pi e/2  < d_i-f_i \leq \pi e/2   \quad \quad \quad  \text{se}  \quad  d_r < f_r   \\  
\end{split}
\]

Dunque i numeri pseudouguali a $f$ sono solo $f$ e $\s f$ (e le loro  traslazioni di $i 2 \pi e$). Per il resto ovviamente laddove $d$ non è pseudomaggiore di $f$  è pseudominore, e viceversa.

Immaginiamo di colorare di blu ogni $d$ pseudomaggiore dell' $f$ fissato, e di rosso ogni $d$ pseudominore.
Osserviamo che l'interpretazione geometrica della 27.3) è che, traslando il piano di $i \pi e$,  deve prodursi un'inversione dei colori in ogni punto.
L'interpretazione geometrica della 27.1) è invece che, ruotando il piano di 180°  attorno a $f$, deve prodursi un'inversione dei colori in ogni punto, a parte due punti speciali: $f$  resta fisso, e $f + i \pi e$ si sovrappone a   $f - i \pi e$, che è lo stesso numero ($\s f$), proprio perché l'asse immaginario è $ i\R /i2\pi e$. Ciò implica che $\s f$ non può essere né pseudomaggiore né pseudominore di $f$. Questa è la seconda ragione per cui già in $E$ avevamo definito   $a \asymp \s a$.

\begin{figure} [h]
\includegraphics[scale=0.5]{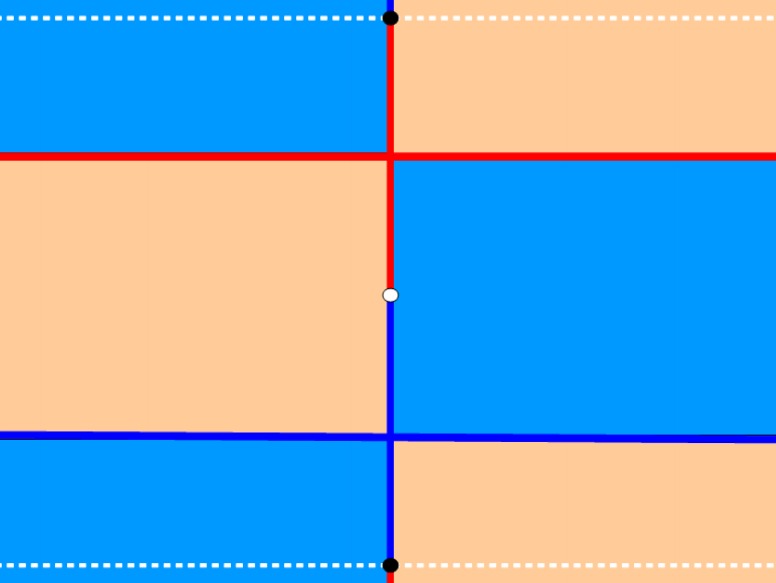}
\centering
\caption{Pseudoordinamento dei complessi. Le linee tratteggiate bianche, perpendicolari all'asse verticale, delimitano una fascia di altezza $ i2 \pi e$.  Tale fascia è divisa in quattro parti uguali da tre linee parallele alle linee bianche: la linea rossa, la linea blu e la linea orizzontale passante per il punto bianco (non disegnata). Rispetto al punto bianco centrale, i punti in blu/azzurro sono pseudomaggiori, i punti in rosso/arancione sono pseudominori e i punti neri sono pseudouguali.}

\end{figure}

A questo punto si possono estendere l'incrementazione e la decrementazione a $\C$ in maniera ovvia.

\[
\begin{split}
29) \quad & d \odot f= \{d+1, \s d +1 \} \quad  d,f \in \C  \quad \text{se}  \quad d \succeq f\\
                 & d \odot f =\{f+1, \s f +1 \}  \quad   d,f \in \C   \quad \text{se} \quad d \prec f\\
\end{split}
\]

Come al solito, il primo valore tra graffe è il valore principale.

\[
\begin{split}
30) \quad & g \oslash  d = g-1 \quad      d,g \in \C       \quad \text{se} \quad d \preceq g-1\\
                 & g \oslash  d = \s (g-1) \quad   d,g \in \C   \quad  \text{se} \quad d \succ g-1\\
\end{split}
\]

Per $d=g-1$ si ha $g \oslash  d = \{h \in \C \mid h \preceq g-1  \}$, nella 30) si è lasciato direttamente il valore principale $g-1$.

La decrementazione è dunque chiusa su $\C$. 
Si erano introdotti i numeri Stigma perché il Sik non era chiuso su $\R$. Abbiamo quindi reso il Sik chiuso su $\R$, poi su $E$ e infine su $\C$. \\ 

\section{Un'interpretazione pratica}
Immaginiamo di avere dei sacchetti chiusi, non trasparenti e molto elastici, che contengono altri sacchetti con le stesse caratteristiche. Aprendo un sacchetto vedremo all'interno un certo numero di altri sacchetti.  Dati $n,m,p \in \N$ associamo a ogni naturale $n$ un sacchetto all'interno del quale sono visibili esattamente $n$ sacchetti (non ci interessa se questi sacchetti contengono a loro volta altri sacchetti).  Definiamo $n$ come "valore" del sacchetto. 
Dati due sacchetti di valori $m$ e $p$ con  $m<p$, diciamo che il sacchetto $m$ è minore del sacchetto $p$. L'incrementazione $m \odot p$ equivale a mettere il sacchetto minore nel maggiore, in questa maniera il valore del maggiore diventerà $p+1$.
La decrementazione $p \oslash m$ equivale invece a togliere un sacchetto $m$ da un sacchetto $p$ che diventerà quindi $p-1$.
Consideriamo che questi sacchetti siano ottenuti chiudendo un fazzoletto rettangolare di lati $a,b,c,d$ con $a,c$ lati corti e $b,d$ lati lunghi. La misura del lato $a$ è uguale per tutti i sacchetti. Il lato $c$ puó essere piatto come gli altri oppure opportunamente dentellato, in modo che guardando il sacchetto chiuso si riconosca se esso è piatto o dentellato.
Associamo ai naturali i sacchetti piatti e agli stigmanaturali quelli dentellati. Possiamo osservare che il sacchetto risultante dall'incrementazione sarà piatto se il risultato è naturale e dentellato se il risultato è stigmanaturale, in quanto quello che resta esterno è il sacchetto pseudomaggiore, in analogia al valore principale dell'incrementazione.
Questo funziona anche con la decrementazione, ad esempio per $ \s m \oslash p = \s (m-1) $  bisogna estrarre dal sacchetto dentellato $\s m$ un sacchetto piatto $p$, in questo modo il primo sacchetto (che resta dentellato) diventerà $ \s (m-1)$.  
L'addizione tra due sacchetti si fa nel modo seguente: si aprono i due sacchetti spianando i fazzoletti rettangolari e lasciandovi sopra i sacchetti contenuti;  se entrambi i fazzoletti sono dentellati si accostano i due lati dentellati, che sono fatti in modo da incastrarsi perfettamente rendendo il nuovo fazzoletto piatto; altrimenti si accostano due lati corti piatti. Poi si cuciono i due lati accostati e si richiude il fazzoletto con dentro tutti i sacchetti. In questo modo il valore del nuovo sacchetto è la somma dei valori precedenti. Inoltre la somma di due sacchetti piatti o di due dentellati è piatta, mentre la somma di un piatto e un dentellato è dentellata, in analogia alla somma di naturali e stigmanaturali.
L'analogia si puó poi estendere anche a sacchetti negativi, reali e immaginari, con opportuni artifici. 

\begin{figure} [h]
\includegraphics[scale=0.5]{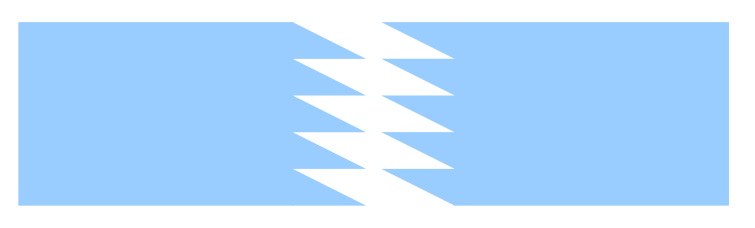}
\centering
\caption{Sacchetti dentellati.}

\end{figure}  

\section *{Ringraziamenti}
Ringrazio sentitamente Rosa Gini e Maurizio Parton per le approfondite discussioni e i loro preziosissimi consigli, Daniele Del Sarto per l'attenta revisione e Gian Franco Romerio per gli utilissimi commenti. \\ \\

\end{document}